\documentclass{article}

\usepackage[utf8x]{inputenc}
\usepackage{amsmath}
\usepackage{amssymb}

\newcommand{\diffop}[1]{\mathrm{d}#1}
\newcommand{\dy}{\diffop{y}}
\newcommand{\df}{\diffop{f}}
\newcommand{\dx}{\diffop{x}}

\newcommand{\dt}{\diffop{t}}

\newcommand{\diff}{\diffop{}}
\newcommand{\equationref}[1]{(\ref{#1})}
\newcommand{\dydx}{\frac{\dy}{\dx}}
\newcommand{\dxdy}{\frac{\dx}{\dy}}
\newcommand{\hdiff}[2]{\mathrm{d}^{#2}#1}

\newcommand{\secondderiv}[2]{\frac{\hdiff{#1}{2}}{\diffop{#2}^2} - \frac{\diffop{#1}}{\diffop{#2}}\frac{\hdiff{#2}{2}}{\diffop{#2}^2}}

\title{Total and Partial Differentials as Algebraically Manipulable Entities}
\author{Maria Isabelle Fite and Jonathan Bartlett}

\begin{document}
\maketitle

\begin{abstract}
    Differential operators usually result in derivatives expressed as a ratio of differentials. For all but the simplest derivatives, these ratios are typically not algebraically manipulable, but must be held together as a unit in order to prevent contradictions. However, this is primarily a notational and conceptual problem. The work of Abraham Robinson has shown that there is nothing contradictory about the concept of an infinitesimal differential operating in isolation. In order to make this system extend to all of calculus, however, some tweaks to standard calculus notation are required. Understanding differentials in this way actually provides a more straightforward understanding of all of calculus for students, and minimizes the number of specialized theorems students need to remember, since all terms can be freely manipulated algebraically.
\end{abstract}


\section{Introduction}
\label{secIntroduction}

Derivatives are usually written in a notation, such as $\dydx$, where the notation implies that there are two distinct values, $\dy$ and $\dx$, at play.
Historically, $\dy$ and $\dx$ were considered infinitesimal values---values so small that they are practically zero, but not quite zero, and often became real numbers when put in ratio with each other.
This understanding was challenged by practitioners who thought that infinitesimal values were insufficiently rigorous to be used in mathematics.

This led to a reconsideration of derivatives using the concept of a limit.
In the limit definition of the derivative, the $\dy$ and $\dx$ terms do not have independent existences, but exist only within the ratio itself.
In this conception, the ratio is merely suggestive of how the derivative was originally produced but does not represent an actual quotient of two distinct values.
The limit definition of the derivative has been reinforced by the fact that treating differentials as distinct values leads to contradictions in many cases.

However, the work of Abraham Robinson in the 1960s showed that there was no fundamental flaw in expanding the number system to include infinitesimals.
The hyperreal numbers are an extension of the real numbers which allows for infinitesimals and infinities to be constructed in a manner equally rigorous with the real numbers.
Additionally, unlike other conceptions of infinities, the hyperreal numbers have an additional advantage that infinitesimals and infinities can be manipulated using arithmetic and algebraic operations.

However, if infinitesimals can be readily considered without contradiction, why does the notation for derivative operations often lead to contradiction?
The flaw here is actually in the notation itself.
Because the notation was not considered factual but merely suggestive, practitioners tended to ignore the problematic cases rather than solve them.
By considering new and more rigorous approaches to notation, a better notation can be developed which  includes infinitesimal values, removes the contradictions, and provides a more straightforward understanding of differential notation and formulas.
In these new formulations, differentials such as $\dy$ and $\dx$ are fully independent, algebraically manipulable entities.

\section{Problem of Separating Differentials in Modern Leibniz Notation}
\label{secProblemOfSeparating}

While the problems that occur when trying to separate differentials in modern Leibniz notation are well-known, it is worth revisiting them briefly.
First of all, it is interesting to note that there are essentially no inconsistencies or contradictions when dealing with first-order total differentials.
For instance, taking the equation $y = x^3$, the derivative is $\dydx = 3x^2$.  
Since the derivative of the inverse function is $\dxdy$, this can be found simply by inverting both sides of the equation, so that $\dxdy = \frac{1}{\dydx} = \frac{1}{3x^2}$.
Likewise, integrating is often preceded by multiplying both sides by a differential, so that $\dydx = 3x^2$ becomes $\dy = 3x^2\,\dx$.

The problems become more apparent on higher-order derivatives.
The typical notation for the second derivative of $y = x^3$ is $\frac{\diffop^2y}{\diffop{x}^2} = 6x$.
However, if the notation were taken seriously, this would be seen as a quotient of the higher-order differential $\diffop^2y$ and the square of $\dx$.
Doing this, however, would break the chain rule.
For instance, if you had $x = t^2$, then you could calculate $\frac{\diffop^2y}{\diffop{t}^2}$ by simply multiplying $\frac{\diffop^2y}{\diffop{x}^2}$ by $\left(\frac{\dx}{\dt}\right)^2$.
Doing so, however, yields an incorrect second derivative of $\frac{\diffop^2y}{\diffop{t}^2} = 24t^4$ rather than the correct $\frac{\diffop^2y}{\diffop{t}^2} = 30t^4$.
This is normally calculated using the chain rule for the second derivative (or higher derivatives using Fa'a di Bruno's formula \cite{tpdiff:johnson2002}).
While the second derivative chain rule works, it provides no algebraic intuition for why it works, and seems to be in conflict with the idea of treating differentials as separable values.

Dealing with partial derivatives brings up innumerable problematic cases even for the first derivative.
If $f$ is a function of $x$ and $y$, and $x$ and $y$ are both functions of $t$, then the total derivative of $f$ with respect to $t$ is $\frac{\df}{\dt} = \frac{\partial f}{\partial x} \frac{\dx}{\dt} + \frac{\partial f}{\partial y} \frac{\dy}{\dt}$.  
Since $x$ is a function of one variable, $\frac{\partial x}{\partial t} = \frac{dx}{dt}$ (likewise for $y$).  
Then the equation becomes $\frac{\df}{\dt} = \frac{\partial f}{\partial x} \frac{\partial x}{\partial t} +  \frac{\partial f}{\partial y} \frac{\partial y}{\partial t}$.
Treating the partial differentials as distinct values, this reduces to $\frac{\df}{\dt} = \frac{\partial f}{\partial t} + \frac{\partial f}{\partial t}$.
\footnote{
A possible objection is that the $∂x$ in $\frac{∂f}{∂x}$ may not be the same infinitesimal as the $∂x$ in $\frac{∂x}{∂t}$.
However, the value of $∂f$ depends on the value of the $∂x$ in $\frac{∂f}{∂x}$, and the value of the $∂x$ in $\frac{∂x}{∂t}$ depends on $∂t$. 
So one could choose the $∂x$s to be equal, and the values of $∂f$ and $∂t$ would adjust accordingly, leaving the values of $\frac{∂f}{∂x}$ and $\frac{∂x}{∂t}$ unchanged.
}
Now that it is expressed in terms of a single variable, $\frac{\df}{\dt} = \frac{\partial f}{\partial t}$, so this yields $\frac{\df}{\dt} = \frac{\df}{\dt} + \frac{\df}{\dt} = 2\frac{\df}{\dt}$.
Dividing both sides by $\frac{\df}{\dt}$ yields the contradiction $1 = 2$.

As will be described, the issues in these problematic cases stem from deficiencies in the notation, not deficiencies in the concept of differentials as infinitesimals nor in the idea that differentials can be considered independently of each other.
By taking a more rigorous approach to the development of the notation of higher order derivatives and partial derivatives, a straightforward notation can be obtained which enables differentials to be considered as fully distinct values.

\section{Historical Formal Definitions of the Derivative}

The derivative of a function measures how the function changes as the independent variable varies. 
For instance, if the derivative of a function $f(x)$ is 3 when $x=5$, that means $f(x)$ is increasing at a rate of 3 units up to every 1 unit across whenever \emph{x} is 5. 
Another way to say the same information is that the function's slope at $x=5$ is $3/1=3$.

Normally, slope is defined with reference to two points. 
When measuring velocity, for instance, which is the ratio of the change in position to the change in time, one would measure two different times with their positions and compare them. 
The derivative attempts to calculate the slope using only one point together with an equation. 
Since only one point is used, the change in $x$ is infinitely small, and so is the change in $y$. 
Different ways of dealing with these infinities lead to different formal definitions of the derivative.

\subsection{Newton's Definition}

Isaac Newton provided one of the first definitions of a derivative in his book \emph{Methodus fluxionum et serierum infinitarum}, or ``The Method of Fluxions and Infinite Series'' in English \cite{tpdiff:newton1736,tpdiff:bell2022}. 
Newton thought of his graphs as being drawn over time, with the $x$-coordinate increasing at a constant speed while the rate of increase in the $y$-coordinate varied. 
A variable's rate of change with respect to time (what we would now call a derivative with respect to time) was called a ``fluxion,'' which was denoted by applying a dot above a variable, such as $\dot x$ (which represents the derivative of $x$ with respect to time) \cite{tpdiff:bell2022}.

To avoid having to define an infinitely small quantity, Newton worked with full derivatives, ratios of infinitesimals. Since Newton assumed all his variables depended on time, he could then switch out the infinitesimal change in $x$ and change in $y$ for the change in $x$ over time and the change in $y$ over time, which were both real numbers. The ratio remained the same, and the infinities were avoided \cite{tpdiff:bell2022}.

\subsection{Leibniz's Definition}

Unlike Newton, Gottfried Leibniz preferred to consider the change in $x$ and the change in $y$ separately. 
He used the notation $\dx$ for an infinitesimal difference in $x$ and $\dy/\dx$ for a ratio of infinitesimals, which represented the slope of a curve at a point. 
Leibniz considered d an operator, with $\dx = \diffop(x)$ being the output of d acting on the variable $x$. 
This allowed him to apply d more than once, resulting in $\mathrm{d}^2x = \diffop(\diffop(x))$, $\mathrm{d}^3x = \diffop(\diffop(\diffop(x)))$, and so on. 
Just like $\dx$ was infinitely smaller than $x$, Leibniz said $\mathrm{d}^nx$ was infinitely smaller than $\mathrm{d}^{n-1}x$ \cite{tpdiff:bell2022}.

Although his calculus relied on the concept of an infinitesimal, Leibniz regarded infinitesimals as only ``purely ideal entities\ldots useful fictions, introduced to shorten arguments and aid insight'' \cite{tpdiff:bell2022}.
However, Leibniz was never able to rigorously define his infinitesimals nor how they behaved.
Therefore, while they seemed to work well, the lack of clarity caused some skeptics to regarded them with suspicion, ridiculing them as ``ghosts of departed quantities'' \cite{tpdiff:berkeley}.

\subsection{Delta-Epsilon (Limit) Definition}

Concerns about the fishy nature of infinitesimals, treated like nonzero numbers when dividing but also like zero when adding, led to the reformulation of calculus using the idea of limits. 
The limit of $f(x)$ as $x$ approaches $a$ is the value $f(x)$ approaches as $x$ becomes closer to $a$.

More precisely, the limit of $f(x)$ as $x$ approaches $a$ is $L$ if for any given positive number $\epsilon$ there is a corresponding
positive number $\delta$ such that the difference between $f(x)$ and $L$ is less than $\epsilon$ whenever the difference between $x$ and $a$ is less than $\delta$ \cite{tpdiff:briggs2019}.

Limits can then be used to define the derivative of a function $f(x)$ as 
\begin{equation}
f^\prime (x) = \lim\limits_{h \to 0}\frac{f(x + h) - f(x)}{h}
\end{equation}

When limits are used to define a derivative, it makes no sense to pull apart the change in $x$ and the change in $y$, as both the limit of the numerator and the limit of the denominator evaluate to zero, and division by zero is undefined.

\section{Hyperreal Numbers and the Definition of the Derivative}

While the limit definition of a derivative solves the philosophical problems of infinitesimals, it does not allow the change in $y$ to be separated from the change in $x$. This led Abraham Robinson to return to Leibniz's infinitesimals in 1958, putting them on a new set-theoretic foundation and creating the field of nonstandard analysis \cite{tpdiff:bell2022}.

While there are different ways to construct hyperreal numbers, the approach we will take here is based on the set theory approach described by Herrmann in \cite{tpdiff:herrmann2010}, with many of the definitions taken from there as well.
We will begin by describing hyperreal numbers (including infinitesimals), and then describe the differential operator as being an operator that can be applied using infinitesimals.

For defining the infinitesimals, the core idea is to take the set of all infinitely long sequences of real numbers, denoted $ℝ^ℕ$. 
Some of these sequences match other sequences so closely they can be considered equivalent. 
Each real number is then assigned to a set of equivalent sequences. 
Then, some of the remaining sets of equivalent sequences can be assigned to infinitesimals. 
Finally, all the operations normally done on real numbers can be translated to operations between sets of equivalent
sequences.

\subsection{Filters, the Cofinite Filter, and Free Ultrafilters: Defining Big Enough}

A \emph{filter} provides a way to classify subsets of a set as either big enough or not big enough. 

Let $X$ be a nonempty set. 
A nonempty subset $F$ of the set of all subsets of $X$ is a proper filter on $X$ if and only if:
\begin{align} 
&\text{(i) for each } A, B ∈ F,\: A ∩ B ∈ F \\
&\text{(ii) if }A ⊂ B ⊂ X \text{ and } A ∈ F \text{, then } B ∈ F \\
&\text{(iii) }∅ \notin F
\end{align}

The \emph{cofinite filter C} is defined as 
\begin{equation}
C = x \:|\: (x ⊂ X) \text{ and } (X − x) \text{ is finite}
\end{equation} where $X$ is an infinite set.
$C$ is called the cofinite filter because a subset $x$ of $X$ gets to be in the filter $C$ if and only if $X$ without $x$ is a finite set. 
$C$ gives a mathematical way to define whether an infinite set is considered big enough. 

For instance, if $C$ is the cofinite filter on ℝ, the real numbers, the set of all integers ℤ is not big enough to be in $C$, even though it is an infinite subset of ℝ, because there are infinitely many real numbers that are not integers. 
However, ℝ*, the real numbers excluding zero, is big enough to be a member of $C$, because there is only one real number, zero, that is not in the real numbers excluding zero.

An \emph{ultrafilter} is the biggest filter on a given infinite set $X$.
An ultrafilter that has $C$ as a subset is called a \emph{free ultrafilter}.

\subsection{Equivalence Classes of $ℝ^ℕ$: Classifying Equivalent Sequences Together}

Let $ℝ^ℕ$ represent the set of all sequences with domain ℕ and range values in ℝ.
Let $A$ and $B$ be two sequences in $ℝ^ℕ$.
$A$ is said to be equivalent to $B$ $(A =_U B)$ if a sufficiently large number of their elements match, or
\begin{equation}
A =_U B \iff n \:|\: \{A_n = B_n\} = S ∈ U
\end{equation}
 
The free ultrafilter $U$ determines whether the set of matching elements is big enough.

This relation $=_U$ is an equivalence relation on
$ℝ^ℕ$, so it can partition $ℝ^ℕ$ into equivalence classes. 
Each equivalence class $[A]$ contains all the sequences in $ℝ^ℕ$ that are equivalent to $A$, including $A$ itself.

The set of all these equivalence classes is called the set of the hyperreal numbers, denoted ${}^*ℝ$.

\subsection{Connecting the Real Numbers to the Hyperreals}

We can define a function $f$ that takes each $x ∈ ℝ$ and gives the unique $[R]$, where $\{n \: | \: R_n = x\} ∈ U$. 
This function \emph{f} assigns to each real number $x$ a hyperreal number $[R]$, namely that set of all sequences where a sufficiently large number of each sequence's elements is $x$.
Often, $f(x)$ is represented by ${}^*x$.
For instance, the hyperreal ${}^*3$ is the set of all
sequences equivalent ($=_U$) to $\{3, 3, 3, \ldots\}$.

Most applications of math use real numbers, so it is helpful to define the subset of the hyperreals that corresponds to the real numbers. 
The image of a subset $X$ of ℝ under $f$ is denoted ${}^\sigma X$. Each hyperreal number ${}^*x$ in ${}^\sigma X$ corresponds to a real number $x$ in $X$.
Since $ℝ$ is a subset of $ℝ$, ${}^\sigma ℝ$ is the subset of the hyperreals that corresponds to the real numbers.

\subsection{Operations on the Hyperreals}

In order for algebra in ${}^*ℝ$ to replace algebra in the real numbers, operations like + and ⋅, among others, have to be defined between members of ${}^*ℝ$. 
It is also useful to define the relation ≤ and the absolute value function.

Let $a$, $b$, and $c$ be elements of ${}^*ℝ$, and let ${}^*+: {}^*ℝ → {}^*ℝ$ be defined as 
\begin{equation}
a {}^*+ b = c \iff \{n \: | \: A_n + B_n = C_n \} ∈ U 
\end{equation}
for any $A_n \in a$, $B_n \in b$, and $C_n \in c$.
That is, the sum of 2 elements of ${}^*ℝ$, $a$ and $b$, are equal to another element of ${}^*ℝ$, $c$, if and only if a sufficiently large number of the elements of the sequences $A_n + B_n$ and $C_n$ match, for any sequence $A_n$ in $a$, $B_n$ in $b$, and $C_n$ in $c$. 
Hyperreal multiplication (${}^*⋅$) can be defined similarly.

To construct a hyperreal greater than relation, for each $a = [A], b = [B] ∈ {}^*ℝ$ define 
\begin{equation}
a {}^*≤ b \iff \{n \: | \: A_n ≤ B_n\} ∈ U
\end{equation}

$a {}^*≤ b$ if and only if, given any sequence in $a$ and any sequence in $b$, a sufficiently large number of elements in $a$'s sequence are less than or equal to their corresponding elements in $b$'s sequence.

These operations establish the structure $\langle {}^*ℝ, {}^*+, {}^*⋅, {}^*≤ \rangle$ as a totally ordered field, with $[0]$ as the identity for ${}^*+$ and $[1]$ as the identity for ${}^*⋅$ \cite[pg.~11]{tpdiff:herrmann2010}.

Finally, the absolute value function can be defined for members of $a ∈ {}^*ℝ$ with 
\begin{equation}
    {}^*|a| = |a| = b \iff \{n \: | \: | A_n | = B_n\} ∈ U
\end{equation} 

The absolute value of a hyperreal number $a$ is a hyperreal number $b$ if and only if, given a sequence in $a$ and a sequence in $b$, a sufficiently large number of elements in $b$'s sequence match the absolute value of their corresponding elements in $a$'s sequence.

In summary, +, ⋅, ≤ and the absolute value function, which are defined on the real numbers, can be translated to operations on the hyperreal numbers.

\subsection{Infinitesimals in the Hyperreals}

Not all of the members of ${}^*ℝ$ correspond to real
numbers, because not all sequences of real numbers are constant
sequences. Some of the remaining hyperreals correspond to
infinitesimals.

A hyperreal number $a$ is infinitely large if
\begin{equation}
{}^*x < |a| \text{ for each } {}^*x ∈ {}^\sigma ℝ
\end{equation}
or in other words, if its absolute value is bigger than every hyperreal that corresponds to a real number. 

A hyperreal number $b$ is an infinitesimal or as Newton stated infinitely small if
\begin{equation}
0 ≤ |b| < {}^*x \text{ for each } 0 < x ∈ ℝ.
\end{equation}

Similarly, a hyperreal is an infinitesimal if its absolute value is bigger than or equal to ${}^*0$ and yet smaller than every hyperreal that corresponds to a positive real number.

Notice that ${}^*0$, which is the equivalence class that contains $\{0, 0, 0, \ldots\}$, is the trivial infinitesimal.

For a nontrivial example of an infinitesimal, consider the equivalence class $g$ containing the sequence $\{0, 1, \frac{1}{2}, \frac{1}{3}, \frac{1}{4} \ldots\}$. 
``Then $g ≠ {}^*0$. Now for each $x \in \mathbb{R}^+$ there is some $m \in \mathbb{N}$, $m \neq 0$ such that $0 < \frac{1}{m} < x$.
Thus ${}^*0 < \frac{ {}^*1}{{}^*m} < {}^*x$.
\ldots{} {[}and{]} $g$ is an infinitesimal'' \cite[pg.~17]{tpdiff:herrmann2010}.

\subsection{Division with Infinitesimals}
If infinitesimals are smaller than every real number, can you still divide by them?

Consider a nonzero infinitesimal, say $\epsilon$, and a sequence in $\epsilon$, say $A$. Even if some of $A$'s elements are zeros, $\epsilon \neq {}^*0$, so the set of all zeros in $A$ is not big enough to be in the ultrafilter $U$. 
So, the nonzero elements of $A$ \emph{are} in $U$, since $U$ is an ultrafilter. 
It is then possible to define another sequence $B$ where $B_n = \frac{1}{A_n}$ if $A_n \neq 0$ and $B_n = 0$ if $A_n = 0$. 
$B$ satisfies the property $[A]^* ⋅ [B] = [1]$, and so $[B]$ is the multiplicative inverse of $[A]$. 

In summary, even if there are sequences in $\epsilon$ with zeros, $\frac{[1]}{\epsilon}$ is still defined, and so it is still possible to divide by $\epsilon$
\cite[pg.~11]{tpdiff:herrmann2010}.

\subsection{The Standard and Principal Part Functions}

Hyperreal expressions can be converted into real expressions using the standard part function, $\mathrm{st}()$, which yields the closest real number to the hyperreal expression.
The standard part of an infinitesimal number is always zero.
For infinite values, the standard part yields $+\infty$ or $-\infty$, which is the non-specific infinity indicating that the value is out of range of the real numbers.

The principal part function, $\mathrm{pt}()$, will yield the most significant component of a hyperreal expression \cite{tpdiff:bgn}.
In a hyperreal expression, imagine $\omega$ representing a benchmark infinite value, with $\epsilon = \frac{1}{\omega}$ representing an associated benchmark infinitesimal.
The hyperreal expression $-2\omega^2 + \omega - 5 + 3\epsilon$ represents four different orders of infinity.
The most significant one is $-2\omega^2$, and, thus, it is the principal part.
For the infinitesimal expression $5\epsilon^2 + \epsilon^3$, $5\epsilon^2$ is the principal part.

The principal part of a hyperreal expression is important because non-principal parts, being infinitely less significant than the principal part by definition, do not affect the large-scale behaviors of smooth and continuous functions.

\subsection{Differentials and Derivatives Using Hyperreals}

The derivative of a function $y = f(x)$ using the hyperreals is denoted $\frac{\dy}{\dx}$, the change in $y$ divided by the change in $x$, just like using Leibniz's notation.
However, we can actually define the differentials themselves as infinitesimals, without referring to ratios.

Many have a hard time conceiving of just what a differential is and means.  
It is easy enough to say that a differential is an infinitesimal, but how exactly are individual differentials defined, especially when not being examined in the context of a derivative?  
What exactly does the higher-order notation $\diffop^2y$ mean?

Let us first remember that, in order to be in a relation, two (or more) variables have to be related to each other in some way.
Therefore, we can imagine some variable, let us call it $q$, not explicitly mentioned in the equation, which is in some sense the ``ultimate'' independent variable.

Note that this variable does not need to be explicitly defined.
In fact, it is better if it is not defined explicitly.  
The reason for this is that defining $q$ explicitly means that there is some chance that there exists yet another deeper, more fundamental variable.  
What we are looking for is the deepest, most fundamental, most independent variable.  
Keeping $q$ as a hypothetical independent variable means that our reasoning will continue to hold in the face of finding more and more fundamental quantities.
Our reasoning about an \emph{actual} variable may fail to hold if it is found to not be the fundamental quantity.
We will imagine $q$ to be smoothly increasing by the infinitesimal $\epsilon$.

Since $q$ is the ultimate variable that relates every other variable in the equation, every variable can (theoretically) be written in terms of $q$.
$y$ is actually shorthand for $y(q)$, $x$ is a shorthand for $x(q)$, and so on.
We can then define the differential of an expression (including just a variable) to be the simple difference between the expression at some value $q + \epsilon$ and the expression at some value $q$.
When taking the differential of a variable, we will use the shorthand $\dy$ to mean $\diffop(y)$. 

\begin{equation}
\label{eqDifferentialDefinition}
\dy = \diffop(y) = y(q + \epsilon) - y(q)
\end{equation}

Note that $\dy$ is also a function of $q$ (this fact will become useful when finding the second differential).
Additionally, assuming that $y$ is a smooth and continuous function of $q$, an infinitesimal change in $q$ will lead to an infinitesimal change in in $y$, so $\dy$ will also be infinitesimal.

We can also rearrange (\ref{eqDifferentialDefinition}) and obtain

\begin{equation}
\label{eqDifferentialDefinition2}
y(q + \epsilon) = y(q) + \dy
\end{equation}

These definitions provide a generic definition for the differential and consequent manipulation techniques that can be applied to any expression.
Let us take the simple example $y = x^2$ (which is $y(q) = x(q)^2$)  and apply this differential operator to it.
We will also apply the principal part function at the end in order to simplify the expression to its most consequential portion.

\begin{align*}
y &= x^2 \\
\diffop(y) &= \diffop(x^2) && \text{differential operator}\\
y(q + \epsilon) - y(q) &= x(q + \epsilon)^2 - x(q)^2  && \text{applying \eqref{eqDifferentialDefinition}} \\
\dy &= (x(q) + \dx)^2 - x(q)^2 && \text{applying \eqref{eqDifferentialDefinition2}} \\
\dy &= x(q)^2 + 2x(q)\,\dx + \dx^2 - x(q)^2 && \text{simplifying} \\
\dy &= 2x(q)\,\dx + \dx^2 \\
\dy &= 2x(q)\,\dx && \text{principal part} \\
\dy &= 2x\,\dx && \text{shorthand}
\end{align*}

The second differential is the same process.
It is merely the differential operator applied where differentials are concerned.
$\dy$ is actually $\diffop{(y(q))}$), but we will refer to it as $\dy(q)$ and $\dy(q + \epsilon)$ for a compromise of brevity and clarity.
The notation $\diffop^2y$ will likewise be shorthand for $\diffop(\diffop(y(q)))$.

\begin{align*}
\dy &= 2x\,\dx  \\
\diffop(\dy) &= \diffop(2x\,\dx) && \text{differential operator} \\
 &= 2x(q + \epsilon)\,\dx(q + \epsilon) - 2x(q)\,\dx(q) && \text{applying \eqref{eqDifferentialDefinition}} \\
 &= 2(x(q) + \dx(q))(\dx(q) + \diffop(\dx(q))) - 2x(q)\,\dx(q) && \text{applying \eqref{eqDifferentialDefinition2}} \\
 &= 2x(q)\,\dx(q) + 2x(q)\diffop(\dx(q))  && \text{simplifying} \\ 
 & ~~~+ 2\,\dx(q)^2 + 2\dx\,\diffop(\dx(q)) - 2x(q)\,\dx(q) \\
 &= 2x(q)\diffop(\dx(q)) + 2\dx(q)^2 + 2\dx\,\diffop(\dx(q)) \\
 &= 2x(q)\diffop(\dx(q)) + 2\,\dx(q)^2 && \text{principal part} \\
\diffop^2y &= 2x\,\diffop^2x + 2\,\dx^2 && \text{shorthand}
\end{align*}

This second differential will typically be a second order infinitesimal.
The process can be further repeated for higher order differentials.

The $2x\,\diffop^2x$ term here may be surprising, but the reason for it will become clear in Section~\ref{secextendingmanipulability} when we eliminate the contradictions present in the standard notation for higher-order differentials.

Since all variables in the equation are related to each other, they also share some relationship to $q$.
Therefore, the definition of a differential can be defined universally within an equation without taking into account the specifics of the variables encountered.

Ultimately, taking the differential of a function results in a $\dy$, $\dx$, or some other term.
However, these terms’ definitions are ultimately rooted in this ultimate independent variable $q$, and the results of incrementing it by some hyperreal infinitesimal $\epsilon$.

The derivative, then, is simply a ratio of differentials defined in this way.
While the terminology of ``taking the derivative with respect to $x$'' can still be used, there is no longer anything special about taking the derivative with respect to a variable as opposed to simply dividing by that variable's differential.
Additionally, this expands the ability to take total differentials straightforwardly into multivariable situations, providing that all variables can be, in principle, tied back to some underlying construct like $q$.

\section{Extending the Total Derivative's Algebraic Manipulability}
\label{secextendingmanipulability}

The hyperreal definition of the derivative has several advantages. 
Once hyperreal numbers are defined, the definition of the derivative arises naturally from considering the change in a function when its (theoretical) independent variable changes infinitesimally. 
Unlike the limit definition, the change in $y$ and the change $x$ are separate entities. 
Using hyperreal numbers, we can rigorously define these entities so that they are manipulable using standard algebraic operators. 

However, this requires that we rethink some of the notations from first principles. 
First of all, now that $\dy$ and $\dx$ are reified entities, they now must be considered in applying such rules as the product rule and the quotient rule.
This is straightforward, and the rules are identical to normal calculus rules.
The differential of $x^2\,\dx$ is the result of applying the product rule to the product of $x^2$ and $\dx$, namely $2x\,\dx^2 + x^2\,\diffop^2x$.

When this is taken into account, differentials of any order become algebraically manipulable.

\subsection{The Second Derivative}

Before taking this idea of algebraically manipulable differentials too far, we need to note that the standard notation for the second derivative, $\frac{\diffop^2y}{\diffop{x}^2}$, does not work in this manner.
The problem, here, is that it implies an improper order of operations \cite{tpdiff:bartlett2019}.

Order of operations is very important when doing derivatives.
When doing a derivative, one \emph{first} takes the differential and \emph{then} divides by $\dx$.
The second derivative is the derivative of the first, so the next differential occurs \emph{after the first derivative is complete}, and the process finishes by dividing by $\dx$ again.

However, what does it look like to take the differential of the first derivative?
Basic calculus rules tell us that the quotient rule should be used:
\begin{align*}
\diff\left(\dydx\right) &= \frac{\dx(\diffop(\dy)) - \dy(\diffop(\dx))}{(\dx)^2} \\
                        &= \frac{\hdiff{y}{2}}{\dx} - \frac{\dy}{\dx}\frac{\hdiff{x}{2}}{\dx}
\end{align*}
Then, for the second step, this can be divided by $\dx$, yielding:
\begin{equation}
\label{secondderiv}
\frac{\diffop\left(\frac{\dy}{\dx}\right)}{\dx} = \frac{\hdiff{y}{2}}{\dx^2} - \frac{\dy}{\dx}\frac{\hdiff{x}{2}}{\dx^2}
\end{equation}
This, in fact, yields a notation for the second derivative which is equally algebraically manipulable as the first derivative.
It is not very pretty or compact, but it works algebraically.

The chain rule for the second derivative fits this algebraic notation correctly, provided we replace each instance of the second derivative with its full form (cf.~\equationref{chainrulesecond}):
\begin{equation}
\frac{\hdiff{y}{2}}{\diff{t}^2} - \frac{\diff{y}}{\diff{t}} \frac{\hdiff{t}{2}}{\diff{t}^2} = \left(\secondderiv{y}{x}\right)\left(\frac{\dx}{\dt}\right)^2 + \dydx \left(\secondderiv{x}{t}\right)
\end{equation}
This in fact works out perfectly algebraically.\footnote{Some may be concerned that, in the formula presented in \equationref{secondderiv}, the ratio $\frac{\hdiff{x}{2}}{\dx^2}$ reduces to zero.
However, this is not necessarily true.
The concern is that, since $\frac{\dx}{\dx}$ is always $1$ (i.e., a constant), then $\frac{\hdiff{x}{2}}{\dx^2}$ should be zero.
The problem with this concern is that we are no longer taking $\frac{\hdiff{x}{2}}{\dx^2}$ to be the derivative of $\frac{\dx}{\dx}$.
Using the notation in \equationref{secondderiv}, the derivative of $\frac{\dx}{\dx}$ would be:
\begin{equation}
\label{dxdxderiv}
\frac{\diffop\left(\frac{\dx}{\dx}\right)}{\dx} = \frac{\hdiff{x}{2}}{\dx^2} - \frac{\dx}{\dx}\frac{\hdiff{x}{2}}{\dx^2}
\end{equation}
In this case, since $\frac{\dx}{\dx}$ reduces to $1$, the expression is self-evidently zero.
However, in \equationref{dxdxderiv}, the term $\frac{\hdiff{x}{2}}{\dx^2}$ is not itself necessarily zero, since it is \emph{not} the second derivative of $x$ with respect to $x$.}

\subsection{Higher Order Derivatives}

The notation for the third and higher derivatives can be found using the same techniques as for the second derivative.
To find the third derivative of $y$ with respect to $x$, one starts with the second derivative, takes the differential, and divides by $\dx$:
\begin{equation}
\frac{\diffop\Bigl(\frac{\diffop\left(\dydx\right)}{\dx}\Bigr)}{\dx} = \frac{\diffop\left(\secondderiv{y}{x}\right)}{\dx}  = \frac{\hdiff{y}{3}}{\dx^3} - \frac{\dy}{\dx}\,\frac{\hdiff{x}{3}}{\dx^3} - 3\frac{\hdiff{x}{2}}{\dx^2}\frac{\hdiff{y}{2}}{\dx^2} + 3\frac{\dy}{\dx}\,\frac{(\hdiff{x}{2})^2}{\dx^4}
\end{equation}

Because the expanded notation for the second and higher derivatives is much more verbose than the first derivative, it is often useful for clarity and succinctness to write derivatives using a slight modification of Arbogast's $D$ notation (see \cite{tpdiff:cajori2}) for the total derivative instead of writing it as algebraic differentials.
Here, we will also be subscripting the $D$ with the variable with which the derivative is being taken with respect to and supplying in the superscript the number of derivatives we are taking.  Therefore, where Arbogast would write simply $D$, this notation would be written as $D_x^1$.

Below is the second and third derivative of $y$ with respect to $x$ written using both the enhanced Arbogast notation and as a ratio of differentials.

\begin{align}
D_{x}^2y &= \secondderiv{y}{x} \\
D_{x}^3y &= \frac{\hdiff{y}{3}}{\dx^3} - \frac{\dy}{\dx}\,\frac{\hdiff{x}{3}}{\dx^3} - 3\frac{\hdiff{x}{2}}{\dx^2}\frac{\hdiff{y}{2}}{\dx^2} + 3\frac{\dy}{\dx}\,\frac{(\hdiff{x}{2})^2}{\dx^4}
\end{align}
This gets even more important as the number of derivatives increases.
Each one is more unwieldy than the previous one.
However, each level can be converted to differential notation as follows:
\begin{equation}
D_{x}^{n}y = \frac{\diffop(D_{x}^{n-1}y)}{\dx}
\end{equation}
The advantage of Arbogast's notation over Lagrangian notation are that this modification of Arbogast's notation clearly specifies both the variable/expression whose derivative is being taken and the  variable/expression it is being taken with respect to.

Therefore, when a compact representation of higher order derivatives is needed, this paper will use Arbogast's notation for its clarity and succinctness.  
This notation can be easily expanded to its differentials when necessary for manipulation.

\section{Extending the Partial Derivative's Algebraic Manipulability}

The derivative gives the rate at which a function \emph{f} changes when
$x$ is increased. But what if $f$ depends on both $x$ and
$y$? Imagine a hill where $f$ is the distance above sea level,
$x$ is the distance east from the origin, and $y$ is the
distance north from the origin. To find how $f$ is changing, a
direction to measure the slope must be picked. Along the direction
straight east, only $x$ is changing while $y$ stays constant.
This slope is the partial derivative of $f$ with respect to
$x$, denoted $\frac{\partial f}{\partial x}$, the change in $f$ over the
change in $x$ when $x$ is the only variable allowed to change
\cite[pgs.~940--941]{tpdiff:briggs2019}. A derivative where all the independent
variables are allowed to change is called a total derivative, like the
two-dimensional derivative $\frac{\dy}{\dx}$.
This partial derivative can be formally defined using limits or using hyperreals.

Using limits, the partial derivative of $f(x,y)$ at
the point $(a,b)$ with respect to $x$ is $\lim\limits_{h \to 0}
\frac{f(a + h, b) - f((a,b)}{h}$ \cite[pg.~941]{tpdiff:briggs2019}. Likewise, the partial
derivative of $f(x, y)$ with respect to $x$ is
$\lim\limits_{h \to 0}
\frac{f(x + h, y) - f(x, y)}{h}$. For more than two variables, the
partial derivative of $f(x_1, x_2, \ldots)$ with respect to
$x_1$ is

\begin{equation}
    \frac{∂f}{∂x_1} = \lim_{h \to 0} \frac{f(x_1 + h, x_2, \ldots) - f(x_1, x_2, \ldots)}{h}
\end{equation}

Like the with the total derivative, using limits to define the partial derivative means the change in $f$ and the change in $x$ are not defined separately and must be kept together.
Using hyperreals, the partial derivative of \emph{f} with respect to \emph{x}\textsubscript{1} is 
\begin{equation}
    \frac{∂f}{∂x_1} = \frac{f(x_1 + dx_1, x_2, \ldots) - f(x_1, x_2, \ldots)}{\dx_1}
\end{equation} 

Also, $\dx_1$ can equal $∂x_1$ assuming both of them denote the smallest change in $x_1$ possible. 
This is not an equation in the real numbers; it is an equation in the hyperreals.

Both the numerator and denominator of
$\frac{∂f}{∂x_1}$ have meaning on their own, and they
both are specific hyperreals. So it should be possible to separate the
fraction without problems. 

However, the current notation for $∂f$ does not distinguish between the change in $f$ when
$x_1$ is allowed to change and the change in
$f$ when another variable, say $x_2$, is
allowed to change. In other words, the $∂f$ in
$\frac{∂f}{∂x_1}$ is a different hyperreal from the
$∂f$ in $\frac{∂f}{∂x_2}$, even though they both use the exact same symbol. This can cause problems if the notation is taken seriously (see the contradiction noted in Section~\ref{secProblemOfSeparating}). 
Adding more information to the notation resolves this issue.

The notation for the partial derivative should be changed from $\frac{∂f}{∂x}$ to $\frac{∂(f, x)}{\dx}$ in order to preserve the information in the numerator when the fraction is separated. 

This makes it clear that $∂$ is an operator that takes as an argument not only $f$ but also the choice of which variable to vary. 
The function that $∂$ acts on, in this case $f$, is the first argument of $∂$ and every argument after the first is a variable allowed to change. This can lead to expressions like $∂(f, x, y)$, the change in $f$ when both $x$ and $y$ are allowed to vary.

Using this notation, $\frac{\df}{\dt}$ equals $\frac{∂(f, x)}{\dt} + \frac{∂(f, y)}{\dt}$, not $\frac{\df}{\dt} + \frac{\df}{\dt}$. The contradictions are resolved, and the partial derivative fraction can be separated. 
The numerator and denominator can be moved around just like any other algebraic expression, keeping in mind both of them are hyperreals, so technically any operations on them should be hyperreal operations.

Because the new notation can be algebraically manipulated without contradictions, it makes possible new equations where infinitesimals are not confined to ratios. For instance, the resolved contradiction proof gave the equation $\df = ∂(f, x) + ∂(f, y)$. 
This is reminiscent of one of the conditions for differentiability, $∆f = f_x(a, b) ∆x + f_y(a, b) ∆y + \epsilon_1 ∆x + \epsilon_2 ∆y$, where for fixed $a$ and $b$, $\epsilon_1$ and $\epsilon_2$ are functions that depend only on $∆x$ and $∆y$, with $\epsilon_1, \epsilon_2 → (0, 0)$ as $(∆x,∆y) → (0, 0)$ \cite[pg.~947]{tpdiff:briggs2019}.

Besides simplifying old equations, with the new notation it is possible to consider individual partial changes when building equations, just like considering individual total changes.

The new notation can also denote expressions like $∂(f, x_1, x_2)$, the change in $f(x_1, x_2, x_3)$ when $x_1$ and $x_2$ are allowed to vary, but $x_3$ must stay constant. With the current notation $∂f$, dealing with these situations is clumsy at best.

$∂(f, x_1)$ is an infinitesimal with meaning on its own. It can be defined analogously to Equation \ref{eqDifferentialDefinition}:

\begin{equation}
\label{eqPartialDifferentialDefinition}
    ∂(f, x_1) = f(x_1 + \diffop{x_1}, x_2\ldots) - f(x_1, x_2\ldots)
\end{equation}

The total differential of $f$ is usually defined as the combination of all of the changes in $f$ depending on each variable.
Typically, the total differential of a multivariate function is found using the sum of its partial \emph{derivatives} multiplied by their respective differentials.

\begin{equation}
    \diffop{f}(x_1, x_2 \ldots) = \frac{\partial f}{\partial x_1}\,\diffop{x_1} + \frac{\partial f}{\partial x_2}\,\diffop{x_2} + \ldots  
\end{equation}

Using the new definition of the partial differential, we can rewrite the formula much more straightforwardly, where the total differential is simply a sum of its partial differentials.

\begin{equation}
\label{eqTotalSumOfPartials}
    \diffop{f(x_1, x_2\ldots}) = ∂(f, x_1) + ∂(f, x_2) + \ldots
\end{equation}

\section{Building Differential Formulas}
\label{secswapping}

Using the notation established in this paper, we can build standard calculus formulas in a clear, algebraic manner.
The notation and the formulas will flow directly from the basic truths of calculus and the algebraic reasoning of differentials.

\subsection{The Inverse Function Theorem for Second Derivatives}

The standard inverse function theorem simply states that $\frac{\dx}{\dy} = \frac{1}{\frac{\dy}{\dx}}$.  
In other words, as implied by the algebraic arrangement of its terms, the derivative of $x$ with respect to $y$ is simply the inverse of the derivative of $y$ with respect to $x$.
Using the hyperreal understanding of derivatives allows for a more straightforward way of considering this fact.

More importantly, the new notation for the second derivative likewise allows for a straightforward algebraic construction of an inverse function theorem for the second derivative.
Since the second derivative of $y$ with respect to $x$ is $D^2_x \,y = \frac{\diffop^2 y}{\dx^2} - \frac{\dy}{\dx}\frac{\diffop^2x}{\dx^2}$, then the second derivative of $x$ with respect to $y$ will likewise be $D^2_y \,x = \frac{\diffop^2 x}{\dy^2} - \frac{\dx}{\dy}\frac{\diffop^2y}{\dy^2}$.
Is there a way to construct a formula for converting one to the other?
A simple multiplication by $-\left(\frac{\dx}{\dy}\right)^3$ yields 
\begin{equation*}
-D^2_x\, y \left(\frac{\dx}{\dy}\right)^3 = \frac{\diffop^2x}{\dy^2} - \frac{\diffop^2y}{\dy^2}\frac{\dx}{\dy} 
\end{equation*}

Here, $\frac{\dx}{\dy}$ can be trivially recognized as $\frac{1}{D^1_x y}$, and the right-hand side of the equation can be recognized as $D^2_y x$.
Therefore, this can be rewritten as

\begin{equation}
    -D^2_x\, y \left( \frac{1}{D^1_x y} \right)^3 = D^2_y x
\end{equation}

which is the inverse function theorem for the second derivative.

\subsection{The Chain Rule for the Second Derivative}

The chain rule for the second derivative can also be easily derived from the new notation.
Starting with the notation for the second derivative of $y$ with respect to $x$, we can look at the transformations needed to generate a second derivative of $y$ with respect to $t$.
We will start by multiplying by $\frac{\dx^2}{\dt^2}$ in order to match the leading term to what is needed for the final result.

\begin{align}
    D^2_x\, y &= \frac{\diffop^2 y}{\dx^2} - \frac{\dy}{\dx}\frac{\diffop^2x}{\dx^2} \\
    D^2_x\, y \, (D^1_t\,x)^2 &= \frac{\diffop^2 y}{\dx^2}\frac{\dx^2}{\dt^2} - \frac{\dy}{\dx}\frac{\diffop^2x}{\dx^2}\frac{\dx^2}{\dt^2} \\
     D^2_x\, y \, (D^1_t\,x)^2 &= \frac{\diffop^2 y}{\dt^2} - \frac{\dy}{\dx}\frac{\diffop^2x}{\dt^2} \label{eqCRIntermediate}
\end{align}

In (\ref{eqCRIntermediate}) we see that the leading term is what we want, but the second term is problematic. 
However, it looks a little like the leading term of the second derivative of $x$ with respect to $t$ multiplied by the first derivative of $y$ with respect to $t$. 
Adding that combination to our existing result will yield the desired effect.

\begin{align}
\label{chainrulesecond}
(D^2_x\, y) \, (D^1_t\,x)^2 + (D^1_x\,y)\, (D^2_t\,x) &= \frac{\diffop^2 y}{\dt^2} - \frac{\dy}{\dx}\frac{\diffop^2x}{\dt^2} +  \frac{\dy}{\dx}\frac{\diffop^2x}{\dt^2} - \frac{\dy}{\dx} \frac{\dx}{\dt}\frac{\diffop^2t}{\dt^2} \\
(D^2_x\, y) \, (D^1_t\,x)^2 + (D^1_x\,y)\, (D^2_t\,x) &= \frac{\diffop^2 y}{\dt^2} - \frac{\dy}{\dt} \frac{\diffop^2t}{\dt^2} 
\end{align}

As is evident, the right-hand side is the desired result---the second derivative of $y$ with respect to $t$.

\subsection{The Chain Rule for Multivariate Derivatives}

Building the chain rule for multivariate derivatives is even more straightforward.
Consider a function $f(x,y)$ where $x$ and $y$ are both functions of $t$. 
As noted in (\ref{eqTotalSumOfPartials}), The total change in $f$, $\df$, has two parts: the change due to $x$ changing and the change due to $y$ changing. So,

\begin{equation}
\df = \partial(f,x) + \partial(f, y)
\end{equation}

Dividing both sides by $\dt$,

\begin{equation}
\frac{\df}{\dt} = \frac{\partial(f,x)}{\dt} + \frac{\partial(f, y)}{\dt}
\end{equation}

This is a valid equation, but it is difficult to calculate a value like $\frac{\partial(f,x)}{\dt}$ directly. To make it easier to work with, we can multiply the first term by $\frac{\dx}{\dx}$ and the second by $\frac{\dy}{\dy}$:
\footnote{
Technically, both $\frac{\dx}{\dx}$ and $\frac{\dy}{\dy}$ equal [1], not 1. 
But, since this is an equation in the hyperreals (with hyperreal multiplication), multiplying by the hyperreal multiplication identity doesn't change the value of the right side of the equation.
}

\begin{align}
\frac{\df}{\dt} &= \frac{\partial(f,x)}{\dt} \cdot \frac{\dx}{\dx} + \frac{\partial(f, y)}{\dt} \cdot \frac{\dy}{\dy} \\
&= \frac{\partial(f,x)}{\dx} \cdot \frac{\dx}{\dt} + \frac{\partial(f, y)}{\dy} \cdot \frac{\dy}{\dt}
\end{align}

This is the standard chain rule for multivariate derivatives.

\section{Conclusion}

While treating derivatives as ratios of differentials has been long viewed as problematic, small changes in both the understanding and notation of derivatives straightforwardly leads to algebraically manipulable differentials for both total and partial differentials.
These differentials provide a more straightforward basis for both doing calculus operations and deriving standard calculus rules.
It eliminates exceptions and memorized formulas in favor of simply using algebra with differentials.

Our hope is that the flexibility and freedom of manipulability that this notation allows will both reduce the cognitive load for learning to use differential operators as well as allow for easier exploration of possibilities for practitioners.

\section*{Acknowledgments}

The authors wish to thank Dr. Enrique Valderrama for his comments on early drafts of this manuscript.

\end{document}